\documentclass{article}
\usepackage{amsfonts}
\usepackage{amsmath}

\newtheorem{theorem}{Theorem}

\newtheorem{corollary}[theorem]{Corollary}

\newtheorem{lemma}[theorem]{Lemma}

\newenvironment{proof}[1][Proof]{\noindent\textbf{#1.} }{\ \rule{0.5em}{0.5em}}
\input{tcilatex}
\begin{document}

\title{The adjacency matrix of one type of graph and the Fibonacci numbers}
\author{{\small \ Fatih YILMAZ\thanks{%
e-mail adresses: fyilmaz@selcuk.edu.tr, sbbozkurt@selcuk.edu.tr,
dbozkurt@selcuk.edu.tr}, \c{S}. Burcu BOZKURT, Durmu\c{s} BOZKURT} \\
{\small Selcuk University, Science Faculty Department of Mathematics, 42250 }%
\linebreak \\
{\small Campus} {\small Konya, Turkey}}
\maketitle

\begin{abstract}
Recently there is huge interest in graph theory and intensive study on
computing integer powers of matrices.

In this paper, we investigate relationships between one type of graph and
well-known Fibonacci sequence. In this content, we consider the adjacency
matrix of one type of graph with $2k~(k=1,2,...)$ vertices. It is also known
that for any positive integer $r$, the ($i,j$)th entry of $A^{r}$ ($A$ is
the adjacency matrix of the graph) is just the number of walks from vertex $%
i $ to vertex $j,$ that use exactly $k$ edges.
\end{abstract}


\section{Introduction}

There are many special types of matrices which have great importance\ in
many scientific work. For example matrices of tridiagonal, pentadiagonal and
\linebreak others. These types of matrices frequently appear in
interpolation, \linebreak numerical analysis, solution of boundary value
problems, high order harmonic spectral filtering theory and so on. In \cite%
{4}-\cite{7}, the authors investigated computing integer powers of some type
of these matrices.

Among numerical sequences, the Fibonacci numbers which is defined by the
recurrence $F_{n+2}=F_{n+1}+F_{n}$ for $n\geq 0,$ with initial conditions $%
F_{0}=0$ and $F_{1}=1,$ has achieved a kind of celebrity status. Although
Fibonacci sequence has been studied extensively for hundreds of years, it
remains to fascinating and there always seems to be some amazing properties
aspects that are revealed by looking at it closely $\cite{1}$. Fibonacci
sequence has many applications in diverse fields such as mathematics,
computer science, physics, biology and statistics.

A graph $G=(V,E)$ occurs from two finite sets, with set of vertices
\linebreak $V(G)$ $=\{1,2,\ldots ,n\}$ and set of edges $E(G)=\{e_{1},e_{2},%
\ldots ,e_{m}\}$. Let $G$ be a graph with vertices $v_{1},v_{2},\ldots
,v_{n} $. The adjacency matrix of $G$ is an $n$-square matrix $A$ whose $%
(i,j)$th entry, denoted by $[A]_{i,j},$ is defined by: 
\begin{equation*}
\lbrack A]_{i,j}=\left\{ 
\begin{array}{cl}
1, & \text{if }v_{i}\text{ and }v_{j}\text{ are adjacent} \\ 
0, & \text{otherwise}%
\end{array}%
\right.
\end{equation*}%
Also, it is known that, for any positive integer $r$, the $(i,j)$ entry of $%
A^{r}$ is equal to the number of walks from $v_{i}$ to $v_{j}$ that use
exactly $k$ edges $\cite{3}$.

In $\cite{8},$ the authors consider the number of independent sets in graphs
with two elementary cycles. They described the extremal values of the number
of independent sets using Fibonacci and Lucas numbers.

In $\cite{9},$ the author investigated the relationship between $k$-Lucas
sequence and 1-factors of a bipartite graph.

In $\cite{10},$ the authors consider a new family of $k$-Fibonacci numbers
and \linebreak investigate some properties of the relation and well-known
Fibonacci numbers. They give the well-known Binet's formula 
\begin{equation*}
F_{n}=\dfrac{1}{\sqrt{5}}(\alpha ^{n+1}-\beta ^{n+1})
\end{equation*}
where $\alpha =(1+\sqrt{5})/2$ and $\beta =(1-\sqrt{5})/2$.

In $\cite{11},$ the authors give a generalization for known-sequences and
then they give the graph representations of the sequences. They generalize
Fibonacci, Lucas, Pell and Tribonacci numbers and they show that the
sequences are equal to the total number of $k$-independent sets of special
graphs.

In $\cite{12},$ the author derived an explicit formula which corresponds to
the Fibonacci numbers for the number of spanning trees given below:

\FRAME{ftbphF}{2.3601in}{1.2099in}{0pt}{}{}{tree1.eps}{\special{language
"Scientific Word";type "GRAPHIC";maintain-aspect-ratio TRUE;display
"USEDEF";valid_file "F";width 2.3601in;height 1.2099in;depth
0pt;original-width 11.1976in;original-height 5.7242in;cropleft "0";croptop
"1";cropright "1";cropbottom "0";filename 'tree1.eps';file-properties
"XNPEU";}}

In this paper, we consider one type of graph with $2k$ vertices $(n=2k,$ $%
k=1,2,3,\ldots ,\frac{n}{2})$, each of its connected components of length $2$
and every vertex has degree $1$ or $3$. The components of length $2$ are all
in the form an edge followed by a loop. Clearly:

\qquad \qquad \qquad \FRAME{ihF}{2.4924in}{0.6555in}{0in}{}{}{graf1.eps}{%
\special{language "Scientific Word";type "GRAPHIC";maintain-aspect-ratio
TRUE;display "USEDEF";valid_file "F";width 2.4924in;height 0.6555in;depth
0in;original-width 3.5284in;original-height 0.9072in;cropleft "0";croptop
"1";cropright "1";cropbottom "0";filename 'graf1.eps';file-properties
"XNPEU";}}

\qquad \qquad \qquad \qquad \qquad \qquad \qquad {\small Figure 1}

The adjacency matrix of the graph given in Figure 1 is an $n$-square $(0,1)$%
-block-diagonal matrix whose diagonal blocks have the form $[0,1,1,1],$
which is:%
\begin{equation}
A=\left\{ 
\begin{array}{ll}
a_{i,i+1}=a_{i+1,i}=1, & \text{for }i=1,3,5,\ldots ,n-1 \\ 
a_{i,i}=1, & \text{for }i=2,4,6,\ldots ,n \\ 
0, & \text{otherwise}%
\end{array}%
\right.  \label{15}
\end{equation}

The ($i$, $j$)th entry of $A^{r}$ is just the number of the different paths
from vertex $i$ to vertex $j.$ In other words, the number of the different
paths from vertex $i$ to vertex $j$ corresponds Fibonacci numbers.

\section{Main results}

Let us consider the adjecency matrix of the graph given as in (\ref{15}).
One can \linebreak observe that all integer powers of $A$ are specified to
the famous \linebreak Fibonacci numbers with positive and negative signs.

It is also known that the $r$th ($r\in 
\mathbb{N}
$) power of a matrix is computed by using the known expression $%
A^{r}=TJ^{r}T^{-1}$ \cite{2}, here $J$ is the Jordan form of the matrix and $%
T$ is the transforming matrix. The matrices $J$ and $T$ are obtained using
eigenvalues and eigenvectors of the matrix $A.$ The eigenvalues of $A$ are
the roots of the characteristic equation defined by $\left\vert A-\lambda
I\right\vert =0$ where $I$ is the identity matrix of $n$th order.

Let $P_{n}\left( x\right) $ be the characteristic polynomial of the matrix $%
A $ which is defined in (\ref{15}). Then we can write: 
\begin{equation}
\begin{array}{c}
P_{2}\left( x\right) =x^{2}-x-1 \\ 
P_{4}\left( x\right) =x^{4}-2x^{3}-x^{2}+2x+1 \\ 
P_{6}\left( x\right) =x^{6}-3x^{5}+5x^{3}-3x-1 \\ 
P_{8}\left( x\right) =x^{8}-4x^{7}+2x^{6}+8x^{5}-5x^{4}-8x^{3}+2x^{2}+4x+1
\\ 
\vdots%
\end{array}
\label{33}
\end{equation}

Taking (\ref{33}) into account%
\begin{equation*}
\begin{array}{cc}
P_{n}\left( \lambda \right) =(\lambda ^{2}-\lambda -1)^{\frac{n}{2}} &  \\ 
\text{ \ \ \ \ \ \ \ \ \ \ \ }=[(\lambda -\alpha )(\lambda -\beta )]^{\frac{n%
}{2}} & 
\end{array}%
\end{equation*}%
where $n=2k,$ $(k=1,2,\ldots ),$ $\alpha =\dfrac{1+\sqrt{5}}{2}$\ and $\beta
=\dfrac{1-\sqrt{5}}{2}$. The eigenvalues of the matrix are multiple
according to the order of the matrix $A$. Then the Jordan's form of the
matrix $A$ is:

\begin{eqnarray}
J &=&J_{k}=diag(\underbrace{\alpha ,\ldots ,\alpha ,}\underbrace{\beta
,\ldots ,\beta })  \label{38} \\
&&\ \ \ \ \ \ \ \ \ \ \ \ \ \ \ {\small k\ }\text{times}{\small \ \ \ k}%
\text{ times}  \notag
\end{eqnarray}%
where $k=1,2,3,\ldots ,\frac{n}{2}.$ Let us consider the relation $%
J=T^{-1}AT $ ($AT=TJ$); here $A$ is n$th$ order matrix ($n=2k$, $%
k=1,2,\ldots $), $J$ is the jordan form of the matrix $A$ and $T$ is the
transforming matrix. We will find the transforming matrix $T$. Let us denote
the $j$-th column of $T$ by $T_{j}.$ Then $T=(T_{1},T_{2},\ldots ,T_{n})$ and%
\begin{equation*}
(AT_{1},\ldots ,AT_{n})=(\lambda _{1}T_{1},\ldots ,\lambda _{1}T_{k},\lambda
_{2}T_{k+1},\ldots ,\lambda _{2}T_{2k}).
\end{equation*}

In other words 
\begin{equation}
\begin{array}{ccc}
AT_{1} & = & \lambda _{1}T_{1} \\ 
AT_{2} & = & \lambda _{1}T_{2} \\ 
& \vdots &  \\ 
AT_{k} & = & \lambda _{1}T_{k} \\ 
AT_{k+1} & = & \lambda _{2}T_{k+1} \\ 
AT_{k+2} & = & \lambda _{2}T_{k+2} \\ 
& \vdots &  \\ 
AT_{2k} & = & \lambda _{2}T_{2k}%
\end{array}%
.  \label{20}
\end{equation}

Solving the set of equations system, we obtain the eigenvectors of the
matrix $A:$

\begin{eqnarray}
T &=&\left( \underbrace{%
\begin{array}{cccc}
0 & 0 & \cdots & 1 \\ 
0 & 0 & \cdots & \alpha \\ 
0 & 0 & \cdots & 0 \\ 
0 & 0 & \cdots & 0 \\ 
\vdots & \vdots & {\mathinner{\mkern2mu\raise1pt\hbox{.}\mkern2mu
\raise4pt\hbox{.}\mkern2mu\raise7pt\hbox{.}\mkern1mu}} & \vdots \\ 
0 & 1 & \cdots & 0 \\ 
0 & \alpha & \cdots & 0 \\ 
1 & 0 & \cdots & 0 \\ 
\alpha & 0 & \cdots & 0%
\end{array}%
}\right. \left. \underbrace{%
\begin{array}{ccccc}
0 & 0 & \cdots & 0 & 1 \\ 
0 & 0 & \cdots & 0 & \beta \\ 
0 & 0 & \cdots & 1 & 0 \\ 
0 & 0 & \cdots & \beta & 0 \\ 
\vdots & \vdots & {\mathinner{\mkern2mu\raise1pt\hbox{.}\mkern2mu
\raise4pt\hbox{.}\mkern2mu\raise7pt\hbox{.}\mkern1mu}} & \vdots & \vdots \\ 
0 & 1 & \cdots & 0 & 0 \\ 
0 & \beta & \cdots & 0 & 0 \\ 
1 & 0 & \cdots & 0 & 0 \\ 
\beta & 0 & \cdots & 0 & 0%
\end{array}%
}\right)  \label{45} \\
&&\ \ \ \ \ \ \ \ \ \ \ \ \ k\ \ \ \ \ \ \ \ \ \ \ \ \ \ \ \ \ \ \ \ \ \ \ \
k  \notag
\end{eqnarray}

We shall find the inverse matrix $T^{-1}$ denoting the $i$th row of the
inverse matrix $T^{-1}$ by $T^{-1}=(t_{1},t_{2},\ldots ,t_{n})$ and
implementing the necessary \linebreak transformations, we obtain:%
\begin{equation}
T^{-1}=\frac{1}{\alpha -\beta }\left( 
\begin{array}{ccccccccc}
0 & 0 & 0 & 0 & \cdots & 0 & 0 & -\beta & 1 \\ 
0 & 0 & 0 & 0 & \cdots & -\beta & 1 & 0 & 0 \\ 
\vdots & \vdots & \vdots & \vdots & {\mathinner{\mkern2mu\raise1pt\hbox{.}%
\mkern2mu \raise4pt\hbox{.}\mkern2mu\raise7pt\hbox{.}\mkern1mu}} & \vdots & 
\vdots & \vdots & \vdots \\ 
0 & 0 & -\beta & 1 & \cdots & 0 & 0 & \alpha & -1 \\ 
-\beta & 1 & 0 & 0 & \cdots & \alpha & -1 & 0 & 0 \\ 
0 & 0 & 0 & 0 & \cdots & 0 & 0 & 0 & 0 \\ 
\vdots & \vdots & \vdots & \vdots & {\mathinner{\mkern2mu\raise1pt\hbox{.}%
\mkern2mu \raise4pt\hbox{.}\mkern2mu\raise7pt\hbox{.}\mkern1mu}} & \vdots & 
\vdots & \vdots & \vdots \\ 
0 & 0 & \alpha & -1 & \cdots & 0 & 0 & 0 & 0 \\ 
\alpha & -1 & 0 & 0 & \cdots & 0 & 0 & 0 & 0%
\end{array}%
\right)  \label{55}
\end{equation}%
Using the equalities (\ref{38}), (\ref{45}) and (\ref{55}); we derive the
expression for the $r$th power of the matrix $A:$%
\begin{equation}
A=TJT^{-1}\Rightarrow A^{r}=TJ^{r}T^{-1}=[a_{i,j}(r)]_{n\times n}  \label{75}
\end{equation}%
That is,%
\begin{equation*}
A^{r}=\left\{ 
\begin{array}{l}
a_{i-1,i-1}(r)=\dfrac{1}{\alpha -\beta }(-\beta \alpha ^{r}+\alpha \beta
^{r}) \\ 
a_{i,i}(r)=\dfrac{1}{\alpha -\beta }(\alpha ^{r+1}-\beta ^{r+1}) \\ 
a_{i-1,i}(r)=\dfrac{1}{\alpha -\beta }(\alpha ^{r}-\beta ^{r}) \\ 
a_{i,i-1}(r)=\dfrac{1}{\alpha -\beta }(-\beta \alpha ^{r+1}+\alpha \beta
^{r+1}) \\ 
0,\text{ \ \ \ \ otherwise}%
\end{array}%
\right.
\end{equation*}%
where $i=2,4,6,\ldots ,n$.

\begin{lemma}
Let $A$ be as in (\ref{15}). Then%
\begin{equation*}
\det (A)=(-1)^{k}
\end{equation*}%
where $n=2k,k=1,2,\ldots ,\dfrac{n}{2}.$
\end{lemma}

\begin{proof}
Using Laplace expansion, the determinant can be obtained.
\end{proof}

\begin{corollary}
Let $A=[a_{ij}]$ be $n$-square matrix as in (\ref{15}). Negative integer
\linebreak powers are:%
\begin{equation}
A^{-r}=\left\{ 
\begin{array}{l}
a_{i-1,i-1}(-r)=\dfrac{1}{\alpha -\beta }((-\beta )^{r+1}-(-\alpha )^{r+1})
\\ 
a_{i,i}(-r)=\dfrac{1}{\alpha -\beta }(-\beta (-\alpha )^{r}+\alpha (-\beta
)^{r}) \\ 
a_{i-1,i}(-r)=\dfrac{1}{\alpha -\beta }(-\beta (-\alpha )^{r+1}+\alpha
(-\beta )^{r+1}) \\ 
a_{i,i-1}(-r)=\dfrac{1}{\alpha -\beta }((-\beta )^{r}-(-\alpha )^{r}) \\ 
0,\text{ \ \ \ \ otherwise}%
\end{array}%
\right.
\end{equation}%
where $i=2,4,6,\ldots ,n$ and $r=1,2,\ldots .$
\end{corollary}

\begin{proof}
By (\ref{75}), we can write $J=T^{-1}AT.$\ We also can rewrite 
\begin{equation*}
J^{-1}=(T^{-1}AT)^{-1}=T^{-1}A^{-1}T.
\end{equation*}%
It is known that the jordan matrix for $A^{-1}$ is: 
\begin{eqnarray}
J^{-1} &=&diag(\underbrace{1/\alpha ,\ldots ,1/\alpha ,}\underbrace{1/\beta
,\ldots ,1/\beta })  \label{90} \\
&&\ \ \ \ \ \ \ \ \ \ \QTR{sl}{k\ }\text{times}\ \ \ \ \ \ \ \ \ k\text{
times}  \notag
\end{eqnarray}%
Provided the equalities (\ref{45}), (\ref{55}) and (\ref{90}), the proof can
be seen easily.
\end{proof}

\section{Examples}

We can find the arbitrary integer powers of the matrix $A,$ taking into
account derived expressions. For example, if $k=3$:%
\begin{equation*}
A=\left( 
\begin{array}{cccccc}
0 & 1 & 0 & 0 & 0 & 0 \\ 
1 & 1 & 0 & 0 & 0 & 0 \\ 
0 & 0 & 0 & 1 & 0 & 0 \\ 
0 & 0 & 1 & 1 & 0 & 0 \\ 
0 & 0 & 0 & 0 & 0 & 1 \\ 
0 & 0 & 0 & 0 & 1 & 1%
\end{array}%
\right) .
\end{equation*}%
For $r=4:$%
\begin{equation*}
A^{4}=\left\{ 
\begin{array}{l}
a_{11}(4)=a_{33}(4)=a_{55}(4)=\dfrac{1}{\alpha -\beta }(-\beta \alpha
^{4}+\alpha \beta ^{4})=F_{3} \\ 
a_{22}(4)=a_{44}(4)=a_{66}(4)=\dfrac{1}{\alpha -\beta }(\alpha ^{5}-\beta
^{5})=F_{5} \\ 
a_{12}(4)=a_{34}(4)=a_{56}(4)=\dfrac{1}{\alpha -\beta }(\alpha ^{4}-\beta
^{4})=F_{4} \\ 
a_{21}(4)=a_{43}(4)=a_{65}(4)=\dfrac{1}{\alpha -\beta }(-\beta \alpha
^{5}+\alpha \beta ^{5})=F_{4} \\ 
0,\text{ \ \ \ \ otherwise}%
\end{array}%
\right. \text{ \ \ \ \ \ \ \ \ \ \ \ \ \ \ \ }
\end{equation*}%
For $r=-4$:%
\begin{equation*}
A^{-4}=\left\{ 
\begin{array}{l}
a_{11}(-4)=a_{33}(-4)=a_{55}(-4)=\dfrac{1}{\alpha -\beta }((-\beta
)^{5}-(-\alpha )^{5})=F_{5} \\ 
a_{22}(-4)=a_{44}(-4)=a_{66}(-4)=\dfrac{1}{\alpha -\beta }(-\beta (-\alpha
)^{4}+\alpha (-\beta )^{4})=F_{3} \\ 
a_{12}(-4)=a_{34}(-4)=a_{56}(-4)=\dfrac{1}{\alpha -\beta }(-\beta (-\alpha
)^{5}+\alpha (-\beta )^{5})=-F_{4} \\ 
a_{21}(-4)=a_{43}(-4)=a_{65}(-4)=\dfrac{1}{\alpha -\beta }((-\beta
)^{4}-(-\alpha )^{4})=-F_{4} \\ 
0,\text{ \ \ \ \ otherwise}%
\end{array}%
\right. \text{ \ }
\end{equation*}%
For $r=-5:$%
\begin{equation*}
A^{-5}=\left\{ 
\begin{array}{l}
a_{i-1,i-1}(-5)=\dfrac{1}{\alpha -\beta }((-\beta )^{6}-(-\alpha
)^{6})=-F_{6} \\ 
a_{i,i}(-5)=\dfrac{1}{\alpha -\beta }(-\beta (-\alpha )^{5}+\alpha (-\beta
)^{5})=-F_{4} \\ 
a_{i-1,i}(-5)=\dfrac{1}{\alpha -\beta }(-\beta (-\alpha )^{6}+\alpha (-\beta
)^{6})=F_{5} \\ 
a_{i,i-1}(-5)=\dfrac{1}{\alpha -\beta }((-\beta )^{5}-(-\alpha )^{5})=F_{5}
\\ 
0,\text{ \ \ \ \ otherwise}%
\end{array}%
\right. .\text{ \ \ \ \ \ \ \ \ \ \ \ \ \ \ \ \ \ \ \ \ \ \ \ \ \ \ \ \ \ \
\ \ \ \ \ }
\end{equation*}


\bigskip

\begin{thebibliography}{99}
\bibitem{1} T. Koshy, Fibonacci and Lucas Numbers with Applications,
Wiley-Interscience Publication, 2001.

\bibitem{2} R. Horn, C. Johnson, Matrix Analysis,Cambridge University Press,
\linebreak Cambridge, 1985.

\bibitem{3} J. M$.$ Harris$,$ J. L. Hirst, M. J. Mossinghoff, Combinatorics
and Graph Theory, Second Edition, Springer 2008.

\bibitem{4} J. Rimas, On computing of arbitrary positive integer powers for
one type of symetric pentadiagonal matrices of even order, Appl. Math.
Comput., 203 (2008) 582-591.

\bibitem{5} H. K\i yak, I. Gurses, F. Y\i lmaz, D. Bozkurt, A formula for
computing integer powers for one type of tridiagonal matrix, Hacettepe
Journal of \linebreak Mathematics and Statistics, Volume 39 (3) (2010)
351-363

\bibitem{7} J. Rimas, On computing of arbitrary positive integer powers for
one type of symmetric tridiagonal matrices of even order-I, Appl. Math.
Comput. 168 (2005) 783-787.

\bibitem{8} M. Startek, A. Wloch, I. Wloch, Fibonacci numbers and Lucas
numbers in graphs, Discrete Applied Mathematics, 157 (2009) 864-868.

\bibitem{9} G-Y. Lee, $k$-Lucas numbers and associated bipartite graphs,
Linear Algebra and Its Applications, 320 (2000) 51-61.

\bibitem{10} M. El-Mikkawy, T. Sogabe, A new family of $k$-Fibonacci
numbers, Applied Mathematics and Computation, 215 (2010) 4456-4461.

\bibitem{11} I. Wloch, A. Wloch, Generalized sequences and $k$-independent
sets in graphs, Discrete Applied Mathematics, 158 (2010) 1966-1970.

\bibitem{12} Z. R. Bogdanowicz, Formulas for the Number of Spanning Trees in
a Fan, Applied Mathematical Sciences, Vol. 2, 2008, no:16, 781-786.
\end{thebibliography}
\end{document}